\documentclass{amsart}

\theoremstyle{plain}
\newtheorem{theorem}{Theorem}

\theoremstyle{definition}

\newtheorem{remark}{Remark}

\def\R{{\mathord{I\!\! R}}}
\def\N{\textrm{I\kern-0.21emN}}

\def\RR{{\mathcal R}}
\def \th{\mathrm{th}\,}
\def \ch{\mathrm{ch}\,}
\def \sh{\mathrm{sh}\,}

\def \dd#1#2{\frac{\partial #1}{\partial #2}}

\def \dd2x#1{\frac{\partial ^2 #1}{\partial x^2}}
\def\div{\mbox{{\rm div}}\;}
\def\rot{\mbox{{\rm curl}}\;}
\def\ds{\longrightarrow}
\def\E{{\mathcal E}}

\def\M{{\mathcal M}}
\def\v{\wedge}

\def\dsp{\displaystyle}

\title[Control of travelling walls in a ferromagnetic nanowire ]
{Control of travelling walls in a ferromagnetic nanowire}

\email{}



\author[G. Carbou, S. Labb\'e, E. Tr\'elat]{}

\begin{document}

\maketitle

\centerline{\scshape Gilles Carbou$^1$, St\'ephane Labb\'e$^2$, Emmanuel Tr\'elat$^3$}
{\footnotesize \centerline{}
\centerline{$^1$ MAB, UMR 5466,
CNRS, Universit\'e Bordeaux 1,}
\centerline{ 351, cours de la Lib\'eration, 33405 Talence cedex,
France.}
\centerline{carbou@math.u-bordeaux1.fr}
\centerline{$^2$ Universit\'e Paris-Sud, Labo. Math., Bat. 425, 91405
Orsay Cedex, France}
\centerline{stephane.labbe@math.u-psud.fr}
\centerline{$^3$ Universit\'e d'Orl\'eans,
Math., Labo.\ MAPMO, UMR 6628,}
\centerline{
Route de Chartres, BP 6759, 45067 Orl\'eans Cedex 2}
\centerline{emmanuel.trelat@univ-orleans.fr}
}

\begin{quote}{\normalfont\fontsize{8}{10}\selectfont
{\bfseries Abstract.}We investigate the problem of controlling the magnetic moment in a
ferromagnetic nanowire submitted to an external magnetic field in the
direction of the nanowire. The system is modeled with the one
dimensional Landau-Lifschitz equation. In the absence of control,
there exist particular solutions, which happen to be relevant for
practical issues, called travelling walls. In this paper, we prove
that it is possible to move from a given travelling wall profile to
any other one, by acting on the external magnetic field. The control
laws are simple and explicit, and the resulting trajectories are shown
to be stable.
 
\par}
\end{quote}

\section{Introduction and main result}\label{sec:1}
The most common model used to describe the behavior of ferromagnetic materials, called
micromagnetism, was introduced by W.-F. Brown in the 60's (see \cite{brown}).
It is based on a thermodynamic approach, and the main idea is that equilibrium states of the magnetization minimize a given energy functional, consisting of several components.
The main components, which permit an accurate description of
the behavior of ferromagnetic materials, are the magnetostatic one
(electromagnetism), the exchange one (spin-spin interactions), the anisotropy one
(crystal shape influence) and the Zeeman one (external influences). This point of view permits to
recover the standard dynamical approach of ferromagnetic phenomena, based on the
so-called \textit{Landau-Lifschitz equation}, which was introduced in the 30's in \cite{LL}.
This equation contains a hamiltonian term stemming from the Larmor microscopic spin precession equation, and a purely dissipative term, perpendicular to the precession component and related to the Euler equation of the static energy functional.

More precisely, ferromagnetic materials are characterized by a spontaneous
magnetization described by the magnetic moment $u$ which is a unitary
vector field linking the magnetic induction $B$ with the magnetic
field $H$ by the relation $B=H+u$. The magnetic moment $u$ is solution of the Landau-Lifschitz equation
\begin{equation}
\label{LL3D}
\frac{\partial u}{\partial t}=-u\wedge H_e-u\wedge (u\wedge H_e),
\end{equation}
where the effective field is given by $H_e=\Delta u + h_d(u)+ H_a$,
and  the demagnetizing field $h_d(u)$  is solution of the
magnetostatic equations
\[\div B=\div(H+u)=0 \mbox{ and } \rot H=0,\]
where $H_a$ is an applied magnetic field.
More details on the ferromagnetism model are provided in \cite{brown,Halpern.Labbe:Modelisation,LL,W}.
Existence results have been established for the Landau-Lifschitz equation in
\cite{carbou.fabrie:time,carbou.fabrie:regular,JMR,V}, numerical aspects have been investigated in
\cite{Had,L,Labbe.Bertin:Microwave}, and asymptotic properties have been proved in
\cite{Alouges:Neel,CFG,desimone.kohn.muller.otto:magnetic,riviere.serfaty:compactness,sanchez}.

In this article, we consider an asymptotic one dimensional model of
ferromagnetic nanowire submitted to an applied field along the axis
of the wire. Let $(e_1,e_2,e_3)$ denote the canonical basis of
$\R^3$.  The ferromagnetic nanowire is assumed to have an infinite length, and is represented by
the axis $\R e_1$. The demagnetizing energy is given by
$h_d(u)=-u_2 e_2 - u_3 e_3$ where $u=(u_1,u_2,u_3)$ (see \cite{sanchez} where this formula has been derived using a BKW method, by considering a wire of nonzero diameter, and taking the limit when
the diameter tends to zero). In addition, we assume
that an external magnetic field $\delta(t) e_1$ is
applied along the wire axis. The real-valued function
$\delta(\cdot)$ is our control.

The Landau-Lifschitz equation writes
\begin{equation}
\label{LL1D}
\dsp \frac{\partial u}{\partial t}=-u\wedge h_\delta (u) -u\wedge(u\wedge
h_\delta(u)),
\end{equation}
where
$h_\delta(u)= \frac{\partial^2u}{\partial x^2}-u_2 e_2 - u_3 e_3 + \delta e_1.$

The magnetic field $\delta(\cdot) e_1$ is generated by a coil winding
up around the nanowire.
Note that it would be extremely costly to generate a magnetic field in
other directions, all along the nanowire.
By the way, such other controls do not lead actually
to relevant solutions.

In terms of control system, setting $h(u)=u_{xx}-u_2e_2-u_3e_3$,
Equation (\ref{LL1D}) yields
\begin{equation}\label{contsys}
u_t = -u\wedge h(u) -u\wedge(u\wedge h(u)) - \delta(u\wedge
e_1+u\wedge (u\wedge e_1)).
\end{equation}

For $\delta\equiv 0$, physical experiments demonstrate the existence
of a particular stationary solution, splitting the nanowire into two
parts. The magnetic moment is almost equal to $e_1$ in one of them,
and to $-e_1$ in the other. This fundamental stationary solution,
called a \textit{wall}, is analytically given by
\begin{equation}
M_0(x) =
\left(
\begin{array}{c}
\th{x}\\ 0 \\ \frac{1}{\ch x}\\
\end{array}
\right).
\end{equation}
Here, and throughout the paper, the notations $\ch$, $\sh$, and $\th$,
respectely stand for the hyperbolic cosine, sine, and tangent
functions.

Stability properties of the solution $M_0$ for System (\ref{contsys})
with $\delta\equiv 0$ have been established in \cite{CL}.

When applying a constant magnetic field in the direction $e_1$ (i.e.,
with a constant control function $\delta(\cdot)\equiv\delta$),
physical experiments show a translation/rotation of the above wall
along the nanowire. The corresponding mathematical solution of
(\ref{contsys}), associated with the constant control $\delta$, is
\begin{equation}
\label{Udelta}
u^\delta(t,x)=R_{\delta t}M_0(x+\delta t),
\end{equation}
where
$$R_\theta=\left(
\begin{array}{ccc}
1&0&0\\
\\
0&\cos \theta &-\sin \theta\\
\\
0&\sin\theta & \cos \theta
\end{array}
\right)
$$
is the rotation of angle $\theta$ around the
axis $\R e_1$.
Furthermore, Equation (\ref{contsys}) is invariant with respect to:
\begin{itemize}
\item translations $x\mapsto x-\sigma$, along the nanowire;
\item rotations $R_\theta$ around the axis $e_1$.
\end{itemize}
Hence, denoting $\Lambda=(\theta,\sigma)$, one has a two parameter family
of symmetries given by $M_\Lambda u(t,x)=R_\theta u(t,x-\sigma)$.
Therefore, we have a three-parameters family of particular solutions of
(\ref{contsys}) defined by
\begin{equation}\label{defiwall}
u^{\delta,\theta,\sigma}(t,x)=M_\Lambda u^\delta(t,x)
=R_{\delta t+\theta}M_0(x+\delta t-\sigma)
\end{equation}
and called \textit{travelling wall profiles}.

\begin{theorem}
\label{thm}
There exist $\varepsilon_0>0$ and $\delta_0>0$ such that, for all $\delta_1,\delta_2\in\R$
satisfying $|\delta_i|\leq\delta_0$, $i=1,2$, for all
$\sigma_1,\sigma_2\in\R$, for every $\varepsilon\in(0,\varepsilon_0)$, there exist $T>0$
and a control function
$\delta(\cdot)\in L^\infty(\R,\R)$ such that, for every solution $u$
of (\ref{contsys}) associated with the control $\delta(\cdot)$ and satisfying
\begin{equation}
\exists\theta_1\in\R\ \vert\
\Vert u(0,\cdot)-u^{\delta_1,\theta_1,\sigma_1}(0,\cdot)\Vert_{H^2}
\leq \varepsilon,
\end{equation}
there exists a real number $\theta_2$ such that
\begin{equation}
\Vert u(T,\cdot)-u^{\delta_2,\theta_2,\sigma_2}(T,\cdot)\Vert_{H^2}\leq\varepsilon.
\end{equation}
Moreover,
there exists real numbers $\theta'_2$ and $\sigma_2'$, with
$\vert\theta'_2-\theta_2\vert+\vert\sigma'_2-\sigma_2\vert\leq\varepsilon$, such that
\begin{equation}
\Vert u(t,\cdot)-u^{\delta_2,\theta_2',\sigma_2'}(t,\cdot)\Vert_{H^2}
\underset{t\rightarrow +\infty}{\longrightarrow} 0.
\end{equation}
\end{theorem}

The control law $\delta(\cdot)$ realizing the conclusion of the
theorem is actually given by the piecewise constant function
\begin{equation}\label{defidelta}
\delta(t)=\left\{\begin{array}{ll}
\delta_2-\frac{\sigma_2-\sigma_1}{T} & \textrm{if}\ 0\leq t\leq T,\\
\delta_2 & \textrm{if}\ t\geq T.
\end{array}\right.
\end{equation}
It combines the advantages of being very simple to implement,
and of sharing robustness properties in $H^2$ norm,
as claimed in the theorem.

The time $T$ of the theorem is arbitrary, but must be large enough so that
$$\left\vert\delta_2-\frac{\sigma_2-\sigma_1}{T}\right\vert\leq\delta_0.$$

This theorem shows that the family of travelling wall profiles (\ref{defiwall}) is approximately controllable in $H^2$ norm, locally in $\delta$ and globally in $\sigma$, in time sufficiently large. The controllability property with respect to $\theta$ is not clear. Intuitively the system should not be controllable in $\theta$, however this question is not very relevant from the physical point of view, since it is the position of the wall which is physically interesting. In particular, our result asserts that it is possible to pass approximately (up to the variable $\theta$ from a wall profile $u^{0,\theta,\sigma}$ to any other by means of a scalar control of the form (\ref{defidelta}). This approximate controllability result may have applications for magnetic recording.
Note that, on the one part, an exact controllability result does not seem to be reachable, due to the physical properties of the system, and on the other part, this approximate controllability property is sufficient for practical interest.

Up to now, only the one dimensional case, that is, a ferromagnetic nanowire, has been considered for control applications. What happens in the two dimensional case is an open question.


\section{Proof of Theorem \ref{thm}}
We follow the same lines as in \cite{CL,CL2}, and first express the
Landau-Lifschitz equation in convenient coordinates. This permits to establish
stability properties, and then to derive the result.

\subsection{Expression of the system in adapted coordinates}
The control function (\ref{defidelta}) considered here being piecewise
constant, it suffices to consider Equation (\ref{contsys}) on each
subinterval.
Hence, we assume hereafter that the control function $\delta(\cdot)$ is
constant, equal to $\delta$. Let $u$ be a solution of (\ref{contsys}). Set
$v(t,x)=R_{-\delta t}(u(t,x-\delta t))$. It is not difficult to check that
\begin{equation}\label{LL2}
v_t=-v\v h(v)-v\v(v\v h(v))-\delta(v_x+ v_1 v - e_1).
\end{equation}
Consider the mobile frame $(M_0(x), M_1(x), M_2)$, where $M_1(\cdot)$ and
$M_2$ are defined by
$$ M_1(x)=\begin{pmatrix}
\dsp \frac{1}{\ch x}\\0\\ -\th
x\end{pmatrix}\mbox{ and
}M_2=\begin{pmatrix}0\\1\\0\end{pmatrix}.$$
In what follows, we will prove that $v$ is close to $M_0$. This allows
to decompose $v:\R^+\times \R\ds S^2\subset \R^3$ in the mobile frame
as
$$v(t,x)=\sqrt{1-r_1(t,x)^2-r_2(t,x)^2}M_0(x)+r_1(t,x)M_1(x)+r_2(t,x)M_2 .$$
Easy but lengthy computations show that $v$ is solution of (\ref{LL2})
if and only if $r=\begin{pmatrix}r_1\\r_2\end{pmatrix}$ satisfies
\begin{equation}\label{LL3}
r_t= A r + R_\delta(x,r,r_x,r_{xx}),
\end{equation}
where
\begin{equation}
R_\delta(x,r,r_x,r_{xx}) = -\delta\begin{pmatrix}\ell & 0\\ 0 &
  \ell\end{pmatrix} r +
G(r)r_{xx}+ H_1(x,r)r_x+H_2(r)(r_x,r_x )+ P_\delta(x,r),
\end{equation}
and
\begin{itemize}
\item $A =\begin{pmatrix}
L & L\\ -L & L \end{pmatrix}$
with $L=\partial_{xx} + (1-2\th ^2 x)\mathrm{Id}$;
\item $\ell=\partial_x+\th x\,\mathrm{Id}$;
\item $G(r)$ is the matrix defined by
\[G(r)=\begin{pmatrix}
\dsp \frac{r_1r_2}{\sqrt{1-|r|^2}} & \dsp \frac{r_2^2}{\sqrt{1-|r|^2}}+ \sqrt{1-|r|^2}-1\\
\dsp -\frac{ r_1^2}{\sqrt{1-|r|^2}}-\sqrt{1-|r|^2}+1  & \dsp
-\frac{r_1r_2}{\sqrt{1-|r|^2}}
\end{pmatrix};\]
\item $H_1(x,r)$ is the matrix defined by
\[H_1(x,r)=\frac{2}{\sqrt{1-|r|^2}\,\ch x }
\begin{pmatrix}
r_2\sqrt{1-|r|^2} -r_1r_2^2 & -r_2+r_2r_1^2\\ \\
r_2-r_2^3 & \sqrt{1-|r|^2} r_2 + r_1 r_2^2
\end{pmatrix};
\]
\item $H_2(r)$ is the quadratic form on $\R^2$ defined by
\[H_2(r)(X,X)=
\frac{ (1-|r|^2)X^TX+(r^TX)^2  }
{(1-|r|^2)^{3/2}}
\begin{pmatrix}
\dsp \sqrt{1-|r|^2} r_1 +
r_2\\ \\
\dsp \sqrt{1-|r|^2}
r_2-r_1
\end{pmatrix};\]
\item $ 
P_\delta(x,r) = \begin{pmatrix} P_\delta^1(x,r) \\ \\ P_\delta^2(x,r) \end{pmatrix},
$
with
\begin{equation*}
\begin{split}
P_\delta^1(x,r) = & 2r_2(\sqrt{1-|r|^2}-1)\frac{1}{\ch^2x}
-2r_1r_2\frac{\sh x}{\ch^2x}
-2r_1|r|^2\frac{1}{\ch^2x}\\
& -2r_1^2\sqrt{1-|r|^2}\frac{\sh x}{\ch^2x}
+r_1^3+r_2(1-\sqrt{1-|r|^2})+r_1r_2^2\\
& -\delta\left(\frac{1}{\ch x}(\sqrt{1-|r|^2}-1+r_1^2) + (\sqrt{1-|r|^2}-1)r_1\th x\right)
\end{split}
\end{equation*}
and
\begin{equation*}
\begin{split}
P_\delta^2(x,r) = & -2r_1(\sqrt{1-|r|^2}-1)\frac{1}{\ch^2x}
+2r_1^2\frac{\sh x}{\ch^2x}
-2r_2|r|^2\frac{1}{\ch^2x} \\
& -2r_1r_2\sqrt{1-|r|^2}\frac{\sh x}{\ch^2x}
+r_2|r|^2 \\
& -\delta\left(\frac{1}{\ch x} r_1r_2 + (\sqrt{1-|r|^2}-1)r_2\th x\right).
\end{split}
\end{equation*}
\end{itemize}

It is not difficult to prove that there exists a constant $C>0$ such that,
if $\Vert r\Vert_{\R^2}^2=|r|^2\leq \frac{1}{2}$ and
$\vert\delta\vert\leq 1$, then, for every $x\in\R$, for all $p,q\in\R^2$,
\begin{equation}
\begin{split}
& \Vert R_\delta(x,r,p,q)\Vert_{\R^2}
\\ & \quad  \leq C ( \vert\delta\vert\Vert p\Vert_{\R^2}
 +\Vert r\Vert_{\R^2}^2 \Vert q\Vert_{\R^2}
+ \Vert r\Vert_{\R^2} \Vert p\Vert_{\R^2} + \Vert r\Vert_{\R^2} \Vert p\Vert_{\R^2}^2 
+ \Vert r\Vert_{\R^2}^2 ).
\end{split}
\end{equation}
This a priori estimate shows that $R_\delta(x,r,r_x,r_{xx})$ is a remainder
term in Equation (\ref{LL3}). The rest of the proof relies on
a spectral analysis of the linear operator $A$, so as to establish
stability properties for Equation (\ref{LL3}).

First of all, notice that $L$ is a selfadjoint operator on $L^2(\R)$, of
domain $H^2(\R)$, and that $L=-\ell^*\ell$ with $\ell=\partial_x+\th x\,\mathrm{Id}$
(one has $\ell^*=-\partial_x+\th x\,\mathrm{Id}$).
It follows that $L$ is nonpositive, and that $\ker L=\ker \ell$ is
the one dimensional subspace of $L^2(\R)$ generated by $\frac{1}{\ch
  x}$. In particular, the operator $L$, restricted to the subspace $E=(\ker
  L)^\perp$, is negative.

\begin{remark}\label{remequivnorm}
It is obvious that, on the subspace $E$:
\begin{itemize}
\item the norms $\Vert (-L)^{1/2}f\Vert_{L^2(\R)}$ and $\Vert f\Vert_{H^1(\R)}$
      are equivalent;
\item the norms $\Vert Lf\Vert_{L^2(\R)}$ and $\Vert f\Vert_{H^2(\R)}$
      are equivalent;
\item the norms $\Vert (-L)^{3/2}f\Vert_{L^2(\R)}$ and $\Vert f\Vert_{H^3(\R)}$
      are equivalent.
\end{itemize}
\end{remark}

Writing $A = JL$, with
$$J= \begin{pmatrix} 1 & 1 \\ -1 & 1 \end{pmatrix},$$
it is clear that the kernel of $A$ is $\ker A=\ker L\times\ker L$; it
is the two dimensional space of $L^2(\R^2)$ generated by
$$
a_1(x)=\begin{pmatrix}
0\\ \frac{1}{\ch x}\end{pmatrix}
\quad \textrm{and}\quad
a_2(x)=\begin{pmatrix}
\frac{1}{\ch x}\\ 0\end{pmatrix}.
$$
Moreover, combining the facts that $L_{\vert (\ker
  L)^\perp}$ is negative and that $\mathrm{Spec}\ J=\{1+i,1-i\}$, it
follows that the operator $A$, restricted to the subspace $\E=(\ker A
)^\perp$, is negative.

These facts suggest to decompose
solutions $r$ of (\ref{LL3}) as the sum of an element of $\ker A$
and of an element of $\E$.

To this aim, recall that, since Equation (\ref{LL2}) is invariant with respect to
translations in $x$ and rotations around the axis $e_1$, for every
$\Lambda=(\theta,\sigma)\in\R^2$,
$M_\Lambda(x)=R_\theta M_0(x-\sigma)$ is solution of (\ref{LL2}).
Define
$$
R_\Lambda(x)=
\begin{pmatrix}
 \langle M_\Lambda(x), M_1(x)\rangle
\\
\langle M_\Lambda(x), M_2\rangle
\end{pmatrix},
$$
the coordinates of $M_\Lambda(x)$ in the mobile frame
$(M_1(x),M_2(x))$.

We claim that the mapping
$$
\begin{array}{rcl}
\Psi: \R^2\times \E& \ds & H^2(\R)\\
(\Lambda,W) & \longmapsto & r(x)=R_\Lambda (x)+W(x)
\end{array}
$$
is a diffeomorphism from a neighborhood $\mathcal{U}$ of zero in
$\R^2\times\E$ into a neighborhood $\mathcal{V}$ of zero in $H^2(\R)$.
Indeed, if $r=R_\Lambda+W$ with $W\in \E$, then, by definition,
\begin{equation}\label{ljk}
\langle r,a_1\rangle_{L^2}=\langle R_\Lambda,a_1\rangle_{L^2}
\quad\textrm{and}\quad \langle r,a_2\rangle_{L^2}=
\langle R_\Lambda,a_2\rangle_{L^2}.
\end{equation}
Conversely, if $\Lambda\in\R^2$ satisfies (\eqref{ljk}), then
$W=r-R_\Lambda\in \E$. The mapping $h:\R^2\ds \R^2$, defined by
$h(\Lambda)= \left( \langle R_\Lambda,a_1\rangle_{L^2},
\langle R_\Lambda,a_2\rangle_{L^2} \right) $
is smooth and satisfies $dh(0)=-2\,\mathrm{Id}$, thus is a local
diffeomorphism at $(0,0)$. It follows easily that $\Psi$ is a local
diffeomorphism at zero.

Therefore, every solution $r$ of (\ref{LL3}), as long as it stays\footnote{This
a priori estimate will be a consequence of the stability property derived next.}
in the neighborhood $\mathcal{V}$, can be written as
\begin{equation}
r(t,\cdot)=R_{\Lambda(t)}(\cdot) + W(t,\cdot),
\end{equation}
where $W(t,\cdot)\in \E$ and $\Lambda(t)\in\R^2$, for every $t\geq 0$, and
$(\Lambda(t),W(t,\cdot))\in\mathcal{U}$.
In these new coordinates\footnote{This decomposition is actually quite standard and
has been used e.g.\ in \cite{kapitula:multidimensionnal} to establish stability
properties of static solutions of semilinear parabolic equations, and in
\cite{bertozzi.munch.shearer:stability,roussier:stability} to prove stability
of travelling waves.}, Equation (\ref{LL3}) leads to (see \cite{CL} for
the details of computations)
\begin{equation} \label{sysWLambda}
\begin{split}
W_t(t,x)&=A W(t,x) + \RR(\delta,\Lambda(t),x,W(t,x),W_x(t,x),W_{xx}(t,x)), \\
\Lambda '(t) &=\M(\Lambda(t),W(t,\cdot),W_x(t,\cdot)),
\end{split}
\end{equation}
where
$\RR:\R\times\R^2\times\R\times \left( H^2(\R)\right)^2\times
\left(H^1(\R)\right)^2\times \left( L^2(\R)\right)^2 \ds
\E$ and $\M:\R^2\times \left( H^1(\R)\right)^2 \times \left(
  L^2(\R)\right)^2 \ds \R^2$ are nonlinear
mappings, for which there exist constants $K>0$ and $\eta>0$ such that
\begin{equation}\label{estiR}
\begin{split}
 & \| \RR(\delta,\Lambda,\cdot,W,W_x,W_{xx})\|_{\left(
     H^1(\R)\right)^2} \\
 &\qquad \qquad
 \leq  K \left(  \Vert\Lambda\Vert_{\R^2} + \vert\delta\vert  + \Vert
 W\Vert_{\left( H^2(\R) \right)^2} \right) \Vert W\Vert_{\left(
     H^3(\R)\right)^2} ,
\end{split}
\end{equation}
\begin{equation}\label{estiM}
 | \M(\Lambda,W,W_x)| \leq K \left(  \Vert\Lambda\Vert_{\R^2}
+\Vert W\Vert_{\left(
     H^1(\R)\right)^2}  \right)  \Vert W\Vert_{\left(
     H^1(\R)\right)^2},
\end{equation}
for every $W\in\E$, every $\delta\in\R$, and every $\Lambda\in\R^2$ satisfying
$\Vert\Lambda\Vert_{\R^2}\leq\eta$.

\begin{remark}
Using the fact that $L$ is selfadjoint, it is obvious to prove that
$AW\in\E$, for every $W\in\E$; hence, (\ref{sysWLambda}) makes sense.
\end{remark}


\subsection{Stability properties, and proof of Theorem \ref{thm}}
We are now in position to establish stability properties for system
(\ref{sysWLambda}). Denoting $W=\begin{pmatrix}W_1\\ W_2\end{pmatrix}$,
define on $\left( H^2(\R)\right)^2\times\R^2$ the function
\begin{equation}
\mathcal{V}(W) = \frac{1}{2}\left\Vert \begin{pmatrix} L & 0 \\ 0 &
  L\end{pmatrix} W\right\Vert_{\left( L^2(\R)\right)^2}^2  = 
\frac{1}{2}\Vert LW_1\Vert_{L^2(\R)}^2 +\frac{1}{2} \Vert LW_2\Vert_{L^2(\R)}^2 .
\end{equation}

\begin{remark}\label{remequivnorm2}
It follows from Remark \ref{remequivnorm} that, on the subspace
$\E=(\ker A)^\perp$, $\sqrt{\mathcal{V}(W)}$ is a norm, which is
equivalent to the norm $\Vert W\Vert_{\left( H^2(\R^2) \right) }^2$.
\end{remark}

Consider a solution $(W,\Lambda)$ of (\ref{sysWLambda}), such that
$W(0,\cdot)=W_0(\cdot)$ and $\Lambda(0)=\Lambda_0$.
Since $L$ is selfadjoint, one has
\begin{equation}\label{eq1}
\begin{split}
& \frac{d}{dt} \mathcal{V}(W(t,\cdot)) 
 = \left\langle AW,
\begin{pmatrix} L^2W_1\\ L^2W_2\end{pmatrix}\right\rangle_{\left(
  L^2(\R)\right)^2} \\
& \qquad +
\left\langle \begin{pmatrix} (-L)^{1/2} & 0 \\ 0 & (-L)^{1/2} \end{pmatrix}
\RR(\delta,\Lambda,\cdot,W,W_x,W_{xx})   ,
  \begin{pmatrix} (-L)^{3/2}W_1\\ (-L)^{3/2}W_2\end{pmatrix}\right\rangle_{\left(
    L^2(\R)\right)^2}.
\end{split}
\end{equation}
Concerning the first term of the right-hand side of (\ref{eq1}), one
computes
$$
\left\langle AW,
\begin{pmatrix} L^2W_1\\ L^2W_2\end{pmatrix}\right\rangle_{\left(
  L^2(\R^2)\right)^2}
= - \Vert (-L)^{3/2}W_1\Vert_{\left( L^2(\R)\right)^2}
- \Vert (-L)^{3/2}W_2\Vert_{\left( L^2(\R)\right)^2},
$$
and, using Remark \ref{remequivnorm}, there exists a constant $C_1>0$
such that
\begin{equation}\label{eq2}
\left\langle AW,
\begin{pmatrix} L^2W_1\\ L^2W_2\end{pmatrix}\right\rangle_{\left(
  L^2(\R)\right)^2}
\leq -C_1 \Vert W\Vert_{\left( H^3(\R)\right)^2}^2 .
\end{equation}
Concerning the second term of the right-hand side of (\ref{eq1}), one
deduces from the Cauchy-Schwarz inequality, from Remark
\ref{remequivnorm}, and from the estimate (\ref{estiR}), that
\begin{equation}\label{eq3}
\begin{split}
&
\left\vert
\left\langle \begin{pmatrix} (-L)^{1/2} & 0 \\ 0 & (-L)^{1/2} \end{pmatrix}
\RR(\delta,\Lambda,\cdot,W,W_x,W_{xx})   ,
  \begin{pmatrix} (-L)^{3/2}W_1\\ (-L)^{3/2}W_2\end{pmatrix}\right\rangle_{\left(
    L^2(\R)\right)^2}
\right\vert
\\
& \quad \leq\ 
\Vert \RR(\delta,\Lambda,\cdot,W,W_x,W_{xx}) \Vert_{\left(
    H^1(\R)\right)^2}
\Vert W\Vert_{\left( H^3(\R)\right)^2} \\
& \quad \leq\ 
 K \left(  \Vert\Lambda\Vert_{\R^2} + \vert\delta\vert  + \Vert
 W\Vert_{\left( H^2(\R) \right)^2} \right) \Vert W\Vert_{\left(
     H^3(\R)\right)^2}^2.
\end{split}
\end{equation}
Hence, from (\ref{eq1}), (\ref{eq2}), and (\ref{eq3}), one gets
$$
\frac{d}{dt} \mathcal{V}(W) \leq 
\left( -C_1 +  K \left(  \Vert\Lambda\Vert_{\R^2} + \vert\delta\vert  + \Vert
 W\Vert_{\left( H^2(\R) \right)^2} \right) \right) \Vert W\Vert_{\left(
     H^3(\R)\right)^2}^2.
$$
If the a priori estimate
\begin{equation*}
 \Vert\Lambda(t)\Vert_{\R^2} + \vert\delta\vert  + \Vert
 W(t,\cdot)\Vert_{\left( H^2(\R) \right)^2} \leq \frac{C_1}{2K}
\end{equation*}
holds, then
$$
\frac{d}{dt} \mathcal{V}(W(t,\cdot)) \leq 
 -\frac{C_1}{2} \Vert W(t,\cdot)\Vert_{\left(
     H^3(\R)\right)^2}^2 \leq  -\frac{C_1}{2} \Vert W(t,\cdot)\Vert_{\left(
     H^2(\R)\right)^2}^2 \leq -C_2 \mathcal{V}(W(t,\cdot))
$$
(using Remark \ref{remequivnorm2} for the existence of a constant $C_2>0$). It follows
that there exist constants $C_3>0$ and $C_4>0$ such that, if $\vert\delta\vert + \Vert W(0,\cdot)\Vert_{\left(
     H^2(\R)\right)^2}\leq\frac{C_1}{4K}$ are small enough, and if the a priori estimate
\begin{equation}\label{apriori}
\max_{0\leq s\leq t}\Vert\Lambda(s)\Vert_{\R^2}   \leq \frac{C_1}{4K}
\end{equation}
holds, then
\begin{equation}\label{cesoir}
\Vert W(s,\cdot)\Vert_{\left( H^2(\R) \right)^2} \leq C_3 \mathrm{e}^{-C_4s}\Vert W(0,\cdot)\Vert_{\left(
     H^2(\R)\right)^2},
\end{equation}
for every $s\in[0,T]$, and moreover,
one deduces from  (\ref{sysWLambda}), (\ref{estiM}), and (\ref{cesoir}) that, if the a priori estimate
(\ref{apriori}) holds, then
\begin{equation}\label{bientotfini}
\begin{split}
\Vert\Lambda(t)\Vert_{\R^2} & \leq \Vert\Lambda(0)\Vert_{\R^2} +
\frac{C_1C_3 }{4}\Vert W(0,\cdot)\Vert_{\left(
     H^2(\R)\right)^2}\int_0^t \mathrm{e}^{-C_4s} ds \\
& \qquad\qquad + KC_3^2\Vert W(0,\cdot)\Vert_{\left(
     H^2(\R)\right)^2}^2 \int_0^t \mathrm{e}^{-2C_4s} ds \\
 & \leq \Vert\Lambda(0)\Vert_{\R^2} +
\frac{C_1C_3 }{4C_4}\Vert W(0,\cdot)\Vert_{\left(
     H^2(\R)\right)^2} + K\frac{C_3^2}{2C_4}\Vert W(0,\cdot)\Vert_{\left(
     H^2(\R)\right)^2}^2 .
\end{split}
\end{equation}

From all previous a priori estimates, we conclude that,
if $\vert\delta\vert+\Vert\Lambda(0)\Vert_{\R^2} + \Vert W(0,\cdot)\Vert_{\left(
     H^2(\R)\right)^2}$ is small enough, then $\Vert\Lambda(t)\Vert_{\R^2}$ remains small, for every $t\geq 0$, and $\Vert W(t,\cdot)\Vert_{\left( H^2(\R)\right)^2}$ is exponentially decreasing to $0$.
     
\medskip
     
The first part of the theorem, on the interval $[0,T]$, easily follows from the above considerations.
For the second part, observe that,
from (\ref{sysWLambda}), (\ref{estiM}), and (\ref{cesoir}), one deduces that
$\Vert\Lambda'(t)\Vert_{\R^2}$ is integrable on $[0,+\infty)$, and hence, $\Lambda(t)$ has a limit
in $\R^2$, denoted $\Lambda_\infty=(\theta_\infty,\sigma_\infty)$, as $t$ tends to $+\infty$.
The theorem follows with $\theta_2'=\theta_2+\theta_\infty$ and $\sigma_2'=\sigma_2+\sigma_\infty$.


\end{document}